\documentclass[11pt]{amsart}
\usepackage{palatino}
\usepackage{amssymb}

\theoremstyle{definition}

\newtheorem*{definition*}{Definition}
\newtheorem*{example*}{Example}
\newtheorem*{remark*}{Remark}
\newtheorem*{question*}{Question}
\newtheorem*{problem*}{Problem}
\newtheorem*{note*}{Note}
\newtheorem*{claim*}{Claim}

\newcommand{\R}{\mathbf{R}}

\newcommand{\ra}{\rightarrow}

\renewcommand{\O}{\mathrm{O}}

\begin{document}

\title{The Probability that a Subspace Contains a Positive Vector}
\author{Kent E. Morrison}
\address{California Polytechnic State University, San Luis Obispo, CA 93407}
\curraddr{American Institute of Mathematics, Palo Alto, CA 94306}
\email{kmorriso@calpoly.edu}
\date{April 2, 2010}
\subjclass[2000]{Primary 60D05}
\keywords{Geometric probability, random subspace}
\begin{abstract}
We determine the probability that a random $k$-dimensional subspace of $\R^n$ contains a positive vector.
\end{abstract}

\maketitle 

\large
\renewcommand{\baselinestretch}{1.2}   
\normalsize

For positive integers $k$ and $n$ with $k \leq n$, let $p(n,k)$ denote the probability that a random $k$-dimensional subspace of $\R^n$ contains a positive vector. The aim of this article is to prove
\begin{equation} 
  p(n,k)=  \frac{1}{2^{n-1}}\sum_{j=0}^{k-1} {{n-1} \choose j} .
\end{equation}

First we make the definitions precise. 
A vector $t \in \R^n$ is \emph{positive}  if $t_i \geq 0$ for all $i$ and $t_i >0$ for at least one $i$, and by \emph{random subspace} we mean point in the Grassmann manifold $G(n,k)$ with its natural $O(n)$-invariant probability measure. This measure is constructed by starting with Haar measure on the orthogonal group $\O(n)$, which is bi-invariant and has total mass 1, and then pushing Haar measure down to $G(n,k)$ using the natural projection $\O(n) \twoheadrightarrow G(n,k)$. We also call this Haar measure.

To prove (1) we use a result of J. G. Wendel \cite{Wendel62} showing that $p(n,d)$ is the probability that $n$ random points in $\R^d$ lie in a half-space or, equivalently, that the convex hull of the points does not contain the origin. Let $d=n-k$ be the complementary dimension for our random subspaces. Given points  $z_1,\dots,z_n \in \R^d$  we define the linear map 
\begin{equation*} \hat{z} \colon \R^n \ra \R^d : (t_1,\dots,t_n) \mapsto \sum t_i z_i .
\end{equation*}
Then the convex hull of the $z_i$ contains the origin if and only if $\ker \hat{z}$ contains a positive vector. (The forward implication is immediate. For the converse, suppose $t \in \ker \hat{z}$ is a positive vector. Thus, $\sum t_i z_i=0$. Then $\sum_i (t_i/T) z_i=0$ is a convex combination of the $z_i$, where $T=\sum_i t_i$.)

If the points $z_i$ are random, then with only mild restrictions on their distribution, $\hat{z}$ has maximal rank, and so the kernel of $\hat{z}$ has dimension $k$. This holds, for example, if the $z_i$ are iid with a distribution absolutely continuous with respect to Lebesgue measure. 
But if we further assume that the $z_i$ are drawn from the probability distribution on $\R^d$ for which the components are iid standard normal variables, then $\ker \hat{z}$ is Haar distributed in $G(n,k)$.

To prove this we note that the distribution of $\hat{z}$ is $\O(n)$-invariant and that the kernel map  from the subset of $d \times n$ matrices of maximal rank to the Grassmannian $G(n,k)$ is $\O(n)$-equivariant. In particular, for a $d \times n$ matrix $A$ and an orthogonal matrix $g \in \O(n)$, we have 
$\ker(Ag^{-1})=g(\ker A)$. It follows that the induced probability measure on $G(n,k)$ is $\O(n)$-invariant and must be Haar measure.

Then the probability that  $\ker \hat{z}$ contains a positive vector is the probability that the origin is in the convex hull of the $z_i$, which is $1-p(n,d)$. Finally, $1-p(n,d)=p(n,k)$, which follows from the identity \[ 2^{n-1}=\sum_{j=0}^{n-1} {{n-1} \choose j}. \]
This completes the proof of (1).

(The identity $p(n,k)+p(n,d)=1$ says that almost surely for $V$ in $G(n,k)$ exactly one of the subspaces $V$ and $V^\perp$ contains a positive vector. This is a probabilistic version of the theorem stating that a subspace contains a positive vector if and only if its orthogonal complement does not contain a strictly positive vector, i.e., a vector all of whose components are positive \cite{Roman08}.)


\begin{thebibliography}{9}
\bibitem{Roman08} S. Roman, \emph{Advanced Linear Algebra}, 3rd ed., Springer, New York, 2008.
\bibitem{Wendel62} J. G. Wendel, A problem in geometric probability, \emph{Math. Scand.} 11 (1962) 109--111.

\end{thebibliography}
\end{document}